# SOBOLEV TESTS OF GOODNESS OF FIT OF DISTRIBUTIONS ON COMPACT RIEMANNIAN MANIFOLDS

By P. E. Jupp

*University of St. Andrews*

Classes of coordinate-invariant omnibus goodness-of-fit tests on compact Riemannian manifolds are proposed. The tests are based on Giné's Sobolev tests of uniformity. A condition for consistency is given. The tests are illustrated by an example on the rotation group $SO(3)$.

**1. Introduction.** Although many tests of goodness of fit are available for distributions on the circle, comparatively little work has been done on developing general tests of goodness of fit on spheres and other sample spaces used in directional statistics. Goodness-of-fit tests for specific models include score tests for Fisher distributions within the Kent family [11], Bingham distributions within the Fisher–Bingham family [11], and for von Mises–Fisher distributions within the Fisher–Bingham family [13], as well as omnibus tests for Fisher distributions [6] and for Watson distributions [2]. An overview is given in Section 12.3 of [14]. The only general work on goodness-of-fit tests for directional distributions appears to be that of Beran [1] and of Boulerice and Ducharme [3]. Beran introduced Wald-type tests for certain nested exponential models on spheres, whereas Boulerice and Ducharme considered score tests of goodness of fit of distributions on spheres and projective spaces. Neither Beran's tests nor those of Boulerice and Ducharme are consistent against all alternatives.

For continuous distributions on the real line or the circle, the probability integral transform can be used to derive a test of goodness of fit from each test of uniformity. However, if the sample space is a manifold of dimension greater than 1, then there is no unique coordinate-invariant analogue of the probability integral transform, so that it is not obvious how one can obtain tests of goodness of fit from tests of uniformity. The purpose of this paper is to use the machinery of Giné's [7] Sobolev tests of uniformity to obtain coordinate-invariant omnibus tests of goodness of fit on arbitrary compact Riemannian manifolds. This is in the spirit of the adaptations of Sobolev tests of uniformity by Wellner [17] to get two-sample tests and by Jupp and Spurr [9, 10] to get tests of symmetry and tests of independence. For a large class of Sobolev tests of uniformity (those which are consistent against all alternatives), the corresponding tests of goodness of fit are consistent against all alternatives. Section 2 recalls Giné's Sobolev tests of uniformity. In Section 3 Sobolev tests of goodness of fit are introduced and their basic properties are given. A numerical example on the rotation group $SO(3)$ is presented in Section 4.

**2. Sobolev tests of uniformity.**   Let $M$ be a compact Riemannian manifold. The Riemannian metric determines the uniform probability measure $\mu$ on $M$. The intuitive idea of the Sobolev tests of uniformity is to map the manifold $M$ into the Hilbert space $L^2(M, \mu)$ of square-integrable functions on $M$ by a function $\mathbf{t} : M \to L^2(M, \mu)$ such that, if $x$ is uniformly distributed, then the mean of $\mathbf{t}(x)$ is $\mathbf{0}$.

The standard way of constructing such mappings $\mathbf{t}$ is due to Giné [7] and is based on the eigenfunctions of the Laplacian operator on $M$. For $k \geq 1$, let $E_k$ denote the space of eigenfunctions corresponding to the $k$th eigenvalue, and put $d(k) = \dim E_k$. Then there is a well-defined map $\mathbf{t}_k$ of $M$ into $E_k$ given by

$$\mathbf{t}_k(x) = \sum_{i=1}^{d(k)} f_i(x) f_i,$$

where $\{f_i : 1 \leq i \leq d(k)\}$ is any orthonormal basis of $E_k$. If $\{a_1, a_2, \dots\}$ is a sequence of real numbers such that

$$(2.1) \qquad \sum_{k=1}^{\infty} a_k^2 \, d(k) < \infty,$$

then

$$(2.2) \qquad x \mapsto \mathbf{t}(x) = \sum_{k=1}^{\infty} a_k \mathbf{t}_k(x)$$

defines a mapping $\mathbf{t}$ of $M$ into $L^2(M, \mu)$. The resulting Sobolev statistic evaluated on observations $x_1, \dots, x_n$ on $M$ is

$$T_n = \frac{1}{n} \left\| \sum_{i=1}^{n} \mathbf{t}(x_i) \right\|^2$$

$$= \frac{1}{n} \sum_{i=1}^{n} \sum_{j=1}^{n} \langle \mathbf{t}(x_i), \mathbf{t}(x_j) \rangle,$$

where $\langle \cdot, \cdot \rangle$ denotes the inner product on $L^2(M, \mu)$ given by

$$\langle f, g \rangle = \int_M f(x) g(x) \, d\mu(x),$$

the integration being with respect to the uniform probability measure $\mu$ on $M$. The corresponding Sobolev test rejects uniformity for large values of $T_n$.

The main properties of $T_n$ are the following:

(i) It is defined without recourse to a coordinate system.

(ii) It is invariant under isometries of $M$.

(iii) Its large-sample asymptotic distribution under uniformity is that of a weighted sum of independent $\chi^2$ distributions.

(iv) The corresponding test is consistent against all alternatives if and only if $a_k \neq 0$ for all $k$.

Further details can be found in [7]. A brief outline of Sobolev tests on spheres is given in Section 10.8 of [14]. Many well-known tests of uniformity are Sobolev tests.

## 3. Sobolev tests of goodness of fit.

3.1. *Weighted Sobolev statistics.* Let $\mathcal{F} = \{f(\cdot; \boldsymbol{\theta}) : \boldsymbol{\theta} \in \Theta\}$ be a family of probability density functions on $M$, where the parameter space $\Theta$ is a $p$-dimensional manifold. The null hypothesis to be tested is that the probability density function of the distribution generating the data is in $\mathcal{F}$. Let $\hat{\boldsymbol{\theta}}$ denote the estimate of $\boldsymbol{\theta}$ obtained from independent observations $x_1, \ldots, x_n$ by means of an estimating function $\psi : M \times \Theta \to \mathbb{R}^p$, that is, $\hat{\boldsymbol{\theta}}$ is the root (assumed unique) of

$$\sum_{i=1}^{n} \psi(x_i; \boldsymbol{\theta}) = \mathbf{0}.$$

The intuitive idea behind the Sobolev goodness-of-fit tests to be introduced here is that under the null hypothesis $\hat{\boldsymbol{\theta}}$ is close to $\boldsymbol{\theta}$, so that the expectation

$$E_{\boldsymbol{\theta}} \left[ \frac{1}{f(x; \hat{\boldsymbol{\theta}})} \mathbf{t}(x) \right]$$

is near $\mathbf{0}$, and so therefore is its sample analogue

$$\frac{1}{n} \sum_{i=1}^{n} \frac{1}{f(x_i; \hat{\boldsymbol{\theta}})} \mathbf{t}(x_i).$$

The closeness of the latter to $\mathbf{0}$ can be measured by the weighted Sobolev statistic

$$(3.1) \qquad T_{\mathrm{w}} = \frac{1}{n} \left\| \sum_{i=1}^{n} \frac{1}{f(x_i; \hat{\boldsymbol{\theta}})} \mathbf{t}(x_i) \right\|^2.$$

Thus, $T_w$ is obtained by applying a Sobolev test of uniformity not to the empirical distribution but to the weighted empirical distribution in which each observation $x_i$ is weighted by the reciprocal of the value $f(x_i; \hat{\boldsymbol{\theta}})$ of the fitted density at that point. The null hypothesis is rejected for large values of $T_w$. Significance can be assessed using Monte Carlo simulation from the fitted distribution.

The weighted Sobolev statistic $T_w$ can also be written as

$$T_w = \frac{1}{n} \sum_{i=1}^{n} \sum_{j=1}^{n} \frac{1}{f(x_i; \hat{\boldsymbol{\theta}}) f(x_j; \hat{\boldsymbol{\theta}})} \langle \mathbf{t}(x_i), \mathbf{t}(x_j) \rangle,$$

which is often suitable for computation.

REMARK 1. Any direct sum decomposition $L^2(M, \mu) = E_1 \oplus E_2$ with $E_1$ and $E_2$ orthogonal in $L^2(M, \mu)$ yields a decomposition $\mathbf{t} = \mathbf{t}_1 + \mathbf{t}_2$ with $\mathbf{t}_j(M) \subseteq E_j$ for $j = 1, 2$, and so

$$T_w = T_{w1} + T_{w2},$$

where

$$T_{wj} = \frac{1}{n} \left\| \sum_{i=1}^{n} \frac{1}{f(x_i; \hat{\boldsymbol{\theta}})} \mathbf{t}_j(x_i) \right\|^2 \qquad \text{for } j = 1, 2.$$

Note that $T_{w1}$ and $T_{w2}$ are not necessarily asymptotically independent under the null hypothesis. Any group $G$ of isometries of $M$ gives such a direct sum decomposition $L^2(M, \mu) = E_{M/G} \oplus E_G$ with

$$E_{M/G} = \{ f \in L^2(M, \mu) : f(gx) = f(x), x \in M, g \in G \},$$

$$E_G = \left\{ f \in L^2(M, \mu) : \int_G f(gx) \, d\lambda(g) = 0 \right\},$$

where $\lambda$ is the uniform probability measure on $G$. If the $f(\cdot; \boldsymbol{\theta})$ are invariant under $G$, in that

$$f(gx; \boldsymbol{\theta}) = f(x; \boldsymbol{\theta}) \qquad \text{for } x \in M, \boldsymbol{\theta} \in \Theta, g \in G,$$

then the component $T_{M/G}$ of $T_w$ obtained from $E_{M/G}$ measures the goodness of fit of the data to the corresponding distribution on the quotient space $M/G$, while the component $T_G$ obtained from $E_G$ measures the lack of symmetry under $G$.

REMARK 2. Beran's [1] goodness-of-fit tests on spheres can easily be generalized to general compact Riemannian manifolds as follows. Let $E_1$ and $E_2$ be orthogonal finite-dimensional subspaces of $L^2(M, \mu)$ which are invariant under isometries of $M$. Consider the exponential model with probability density functions of the form

$$(3.2) \qquad f(x; \boldsymbol{\theta}_1, \boldsymbol{\theta}_2) = \exp\{\langle \boldsymbol{\theta}_1, \mathbf{t}_1(x) \rangle + \langle \boldsymbol{\theta}_2, \mathbf{t}_2(x) \rangle - \kappa(\boldsymbol{\theta}_1, \boldsymbol{\theta}_2)\},$$
$$x \in M, \boldsymbol{\theta}_j \in E_j,$$

where $\mathbf{t}_j : M \to E_j$ for $j = 1, 2$ and $\kappa(\boldsymbol{\theta}_1, \boldsymbol{\theta}_2)$ is the normalizing constant. Then Beran's test of goodness of fit of the model obtained by putting $\boldsymbol{\theta}_2 = \mathbf{0}$ in (3.2) rejects this hypothesis for large values of $\langle \hat{\boldsymbol{\theta}}_2, \mathbf{S}_{22.1}{}^{-1} \hat{\boldsymbol{\theta}}_2 \rangle$, where $\hat{\boldsymbol{\theta}}_2$ is a suitable estimate of $\boldsymbol{\theta}_2$ and $\mathbf{S}_{22.1}{}^{-1}$ is the $(2, 2)$-part of the inverse of the sample variance matrix of $(\mathbf{t}_1(x), \mathbf{t}_2(x))$. There is no direct connection between Beran's tests and the Sobolev goodness-of-fit tests introduced here. The large-sample asymptotic distribution of $\langle \hat{\boldsymbol{\theta}}_2, \mathbf{S}_{22.1}{}^{-1} \hat{\boldsymbol{\theta}}_2 \rangle$ is $\chi^2_{\dim E_2}$ and, in contrast to those Sobolev tests of goodness of fit characterized in Theorem 3 below, Beran's tests are not consistent against all alternatives.

Although Boulerice and Ducharme [3] presented their score tests of goodness of fit only for distributions on spheres and projective spaces, the generalization to distributions on general compact Riemannian manifolds is straightforward. Whereas $T_{\mathrm{w}}$ is defined by (3.1), the statistics of Boulerice and Ducharme have the form

$$T_{\mathrm{BD}} = \bar{\mathbf{h}}' \{\mathrm{var}_{\hat{\boldsymbol{\theta}}}(\bar{\mathbf{h}})\}^{-1} \bar{\mathbf{h}},$$

where

$$\bar{\mathbf{h}} = \frac{1}{n} \sum_{i=1}^{n} \left( \frac{1}{\sqrt{f(x_i; \hat{\boldsymbol{\theta}})}} \mathbf{t}(x_i) - E_0[\sqrt{f(x; \hat{\boldsymbol{\theta}})} \mathbf{t}(x)] \right),$$

$E_0[\cdot]$ denoting expectation under the uniform distribution, and only finitely many $a_k$ are nonzero. Thus, whereas $T_{\mathrm{w}}$ is based on a multiplicative transform of $\mathbf{t}(x_i)$ which makes its mean of order $O(n^{-1/2})$ under the null hypothesis, $T_{\mathrm{BD}}$ is based on a standardization of $\mathbf{t}(x_i)$ which makes its mean zero and its variance matrix the identity under the null hypothesis. In contrast to those Sobolev tests of goodness of fit characterized in Theorem 3 below, the tests based on $T_{\mathrm{BD}}$ are not consistent against all alternatives. One way of obtaining such consistency, mentioned on page 159 of [3], is to replace $T_{\mathrm{BD}}$ by

$$T_{\mathrm{BD}}^* = \frac{1}{n} \left\| \sum_{i=1}^{n} \left( \frac{1}{\sqrt{f(x_i; \hat{\boldsymbol{\theta}})}} \mathbf{t}(x_i) - E_0[\sqrt{f(x; \hat{\boldsymbol{\theta}})} \mathbf{t}(x)] \right) \right\|^2,$$

where in (2.2) $a_k \neq 0$ for all $k$. Because of the need to calculate $E_0[\sqrt{f(x; \hat{\boldsymbol{\theta}})} \mathbf{t}(x)]$, $T_{\mathrm{BD}}$ and $T_{\mathrm{BD}}^*$ are more complicated than $T_{\mathrm{w}}$.

3.2. *Large-sample asymptotic properties.* An appropriate setting for large-sample asymptotic results is that in which the mapping $\mathbf{t}$ given by (2.2) is allowed to depend on the sample size $n$. Thus, there is a sequence $\mathbf{t}_{(1)}, \mathbf{t}_{(2)}, \ldots$ of mappings from $M$ into $L^2(M, \mu)$ of the form

$$(3.3) \qquad \mathbf{t}_{(n)}(x) = \sum_{k=1}^{\infty} a_{n,k} \mathbf{t}_k(x),$$

where the sequences $\{a_{n,1}, a_{n,2}, \ldots\}$ of real numbers satisfy

$$(3.4) \qquad \sum_{k=1}^{\infty} (a_{n,k})^2 \, d(k) < \infty.$$

The corresponding goodness-of-fit statistic is the weighted Sobolev statistic (3.1) with $\mathbf{t}$ replaced by $\mathbf{t}_{(n)}$. If $\mathbf{t}_{(1)}, \mathbf{t}_{(2)}, \ldots$ converges to some limit $\mathbf{t}$, then $T_{\mathrm{w}}$ has a limiting distribution. This is made precise in Theorems 1 and 2 below.

Suppose that $x_1, \ldots, x_n$ are independent observations from some distribution $\nu$ on $M$. Let $\hat{\boldsymbol{\theta}}_\nu$ be the value of $\boldsymbol{\theta}$ (assumed unique) such that

$$E_\nu[\psi(x; \boldsymbol{\theta})] = \mathbf{0}.$$

Then, under standard regularity assumptions (e.g., multivariate versions of those in Sections 4.2.2 and 7.2.2 of [16]) the following distributional result holds.

THEOREM 1 (Asymptotic distribution). *Let $\mathbf{t}_{(1)}, \mathbf{t}_{(2)}, \ldots$ and $\mathbf{t}$ be mappings from $M$ into $L^2(M, \mu)$ given by (3.3) and (2.2), corresponding to sequences which satisfy (3.4) and (2.1). If*

$$(3.5) \qquad \sum_{k=1}^{\infty} (a_{n,k} - a_k)^2 d(k) \to 0 \qquad \text{as } n \to \infty,$$

*then*

$$\frac{1}{\sqrt{n}} \sum_{i=1}^{n} \frac{1}{f(x_i; \hat{\boldsymbol{\theta}})} (\mathbf{t}_{(n)}(x_i) - \boldsymbol{\tau}) \xrightarrow{d} N(\mathbf{0}, \boldsymbol{\Sigma}) \qquad \text{as } n \to \infty,$$

*where $\xrightarrow{d}$ denotes convergence in distribution and*

$$\boldsymbol{\Sigma} = \mathrm{var}_\nu \left( \frac{1}{f(x; \hat{\boldsymbol{\theta}}_\nu)} (\mathbf{t}(x) - \boldsymbol{\tau}) - \left( E_\nu \left[ -\frac{\partial \psi(x; \boldsymbol{\theta})}{\partial \boldsymbol{\theta}} \Big|_{\boldsymbol{\theta} = \hat{\boldsymbol{\theta}}_\nu} \right]^{-1} \psi(x; \hat{\boldsymbol{\theta}}_\nu) \right)' \boldsymbol{v} \right)$$

*with*

$$(3.6) \qquad \begin{aligned} \boldsymbol{\tau} &= E_\nu \left[ \frac{1}{f(x; \hat{\boldsymbol{\theta}}_\nu)} \mathbf{t}(x) \right], \\ \boldsymbol{v} &= E_\nu \left[ \frac{1}{f(x; \hat{\boldsymbol{\theta}}_\nu)} \frac{\partial l(\boldsymbol{\theta}; x)}{\partial \boldsymbol{\theta}'} \Big|_{\boldsymbol{\theta} = \hat{\boldsymbol{\theta}}_\nu} (\mathbf{t}(x) - \boldsymbol{\tau}) \right], \end{aligned}$$

$l(\boldsymbol{\theta}; x)$ *denoting the log likelihood of* $\boldsymbol{\theta}$ *based on a single observation* $x$.

PROOF. Taylor expansion of $\sum_{i=1}^{n} \psi(x_i; \boldsymbol{\theta})'$ about $\hat{\boldsymbol{\theta}}_\nu$ gives

$$\sqrt{n}(\hat{\boldsymbol{\theta}} - \hat{\boldsymbol{\theta}}_\nu) = k_\nu(\hat{\boldsymbol{\theta}}_\nu)^{-1} \frac{1}{\sqrt{n}} \sum_{i=1}^{n} \psi(x_i; \hat{\boldsymbol{\theta}}_\nu)' + O_P(n^{-1/2}),$$

where $\psi(x_i; \boldsymbol{\theta})$ is regarded as a row vector and

$$k_\nu(\hat{\boldsymbol{\theta}}_\nu) = E_\nu \left[ -\frac{\partial \psi(x; \boldsymbol{\theta})'}{\partial \boldsymbol{\theta}} \Big|_{\boldsymbol{\theta} = \hat{\boldsymbol{\theta}}_\nu} \right].$$

Then

$$\frac{1}{\sqrt{n}} \sum_{i=1}^{n} \frac{1}{f(x_i; \hat{\boldsymbol{\theta}})} (\mathbf{t}_{(n)}(x_i) - \boldsymbol{\tau})$$

$$= \frac{1}{\sqrt{n}} \sum_{i=1}^{n} \frac{1}{f(x_i; \hat{\boldsymbol{\theta}}_\nu)} (\mathbf{t}_{(n)}(x_i) - \boldsymbol{\tau})$$

$$+ \frac{1}{\sqrt{n}} \sum_{i=1}^{n} \left( \frac{1}{f(x_i; \hat{\boldsymbol{\theta}})} - \frac{1}{f(x_i; \hat{\boldsymbol{\theta}}_\nu)} \right) (\mathbf{t}_{(n)}(x_i) - \boldsymbol{\tau})$$

$$= \frac{1}{\sqrt{n}} \sum_{i=1}^{n} \frac{1}{f(x_i; \hat{\boldsymbol{\theta}}_\nu)} (\mathbf{t}_{(n)}(x_i) - \boldsymbol{\tau})$$

$$- (\hat{\boldsymbol{\theta}} - \hat{\boldsymbol{\theta}}_\nu)' \frac{1}{\sqrt{n}} \sum_{i=1}^{n} \left\{ \frac{1}{f(x_i; \hat{\boldsymbol{\theta}}_\nu)} \frac{\partial l(\boldsymbol{\theta}; x_i)}{\partial \boldsymbol{\theta}'} \Big|_{\boldsymbol{\theta} = \hat{\boldsymbol{\theta}}_\nu} (\mathbf{t}_{(n)}(x_i) - \boldsymbol{\tau}) - \boldsymbol{v} \right\}$$

$$- \sqrt{n}(\hat{\boldsymbol{\theta}} - \hat{\boldsymbol{\theta}}_\nu)' \boldsymbol{v} + O_P(n^{-1/2})$$

$$= \frac{1}{\sqrt{n}} \sum_{i=1}^{n} \frac{1}{f(x_i; \hat{\boldsymbol{\theta}}_\nu)} (\mathbf{t}_{(n)}(x_i) - \boldsymbol{\tau}) - \sqrt{n}(\hat{\boldsymbol{\theta}} - \hat{\boldsymbol{\theta}}_\nu)' \boldsymbol{v} + O_P(n^{-1/2})$$

$$= \frac{1}{\sqrt{n}} \left\{ \sum_{i=1}^{n} \frac{1}{f(x_i; \hat{\boldsymbol{\theta}}_\nu)} (\mathbf{t}_{(n)}(x_i) - \boldsymbol{\tau}) - \sum_{i=1}^{n} \psi(x_i; \hat{\boldsymbol{\theta}}_\nu) (k_\nu(\hat{\boldsymbol{\theta}}_\nu)^{-1})' \boldsymbol{v} \right\}$$

$$+ O_P(n^{-1/2}).$$

Since $\mathbf{t}$ and $\psi$ are continuous and $M$ is compact, application of the Hilbert space version of the limit theorem for triangular arrays (for the univariate version, see, e.g., Section 1.9.3 of [16]) to

$$\frac{1}{\sqrt{n}} \sum_{i=1}^{n} \left\{ \frac{1}{f(x_i; \hat{\boldsymbol{\theta}}_\nu)} (\mathbf{t}_{(n)}(x_i) - \boldsymbol{\tau}) - \psi(x_i; \hat{\boldsymbol{\theta}}_\nu) (k_\nu(\hat{\boldsymbol{\theta}}_\nu)^{-1})' \boldsymbol{v} \right\}$$

shows that, as $n \to \infty$,

$$\frac{1}{\sqrt{n}} \sum_{i=1}^{n} \frac{1}{f(x_i; \hat{\boldsymbol{\theta}})} (\mathbf{t}_{(n)}(x_i) - \boldsymbol{\tau}) \xrightarrow{d} N(\mathbf{0}, \boldsymbol{\Sigma}),$$

where

$$\boldsymbol{\Sigma} = \mathrm{var}_{\nu} \left( \frac{1}{f(x; \hat{\boldsymbol{\theta}}_{\nu})} (\mathbf{t}(x) - \boldsymbol{\tau}) - \psi(x; \hat{\boldsymbol{\theta}}_{\nu})(k_{\nu}(\hat{\boldsymbol{\theta}}_{\nu})^{-1})' \boldsymbol{\upsilon} \right). \qquad \square$$

The next two results are straightforward consequences of Theorem 1.

THEOREM 2 (Asymptotic null distribution).    *Under the null hypothesis, if* (3.5) *holds, then:*

(i)  $\boldsymbol{\tau} = \mathbf{0}$, *where* $\boldsymbol{\tau}$ *is defined by* (3.6).

(ii)  *The distribution of* $T_{\mathrm{w}}$ *tends as* $n \to \infty$ *to that of* $\|\mathbf{Z}\|^2$, *where* $\mathbf{Z}$ *is a random element of* $L^2(M, \mu)$ *with* $\mathbf{Z} \sim N(\mathbf{0}, \boldsymbol{\Sigma}_0)$ *and*

$$\boldsymbol{\Sigma}_0 = \mathrm{var}_{\nu} \left( \frac{1}{f(x; \hat{\boldsymbol{\theta}}_{\nu})} \mathbf{t}(x) - \psi(x; \hat{\boldsymbol{\theta}}_{\nu}) \left\{ E_{\nu} \left[ -\frac{\partial \psi(x; \boldsymbol{\theta})'}{\partial \boldsymbol{\theta}} \Big|_{\boldsymbol{\theta} = \hat{\boldsymbol{\theta}}_{\nu}} \right]^{-1} \right\}' \boldsymbol{\upsilon} \right)$$

*with*

$$\boldsymbol{\upsilon} = E_0 \left[ \frac{\partial l(\boldsymbol{\theta}; x)}{\partial \boldsymbol{\theta}'} \Big|_{\boldsymbol{\theta} = \hat{\boldsymbol{\theta}}_{\nu}} \mathbf{t}(x) \right],$$

$E_0[\cdot]$ *denoting expectation under the uniform distribution.*

In general, even for quite simple models, the matrices $\boldsymbol{\Sigma}$ and $\boldsymbol{\Sigma}_0$ in Theorems 1 and 2 do not admit simple explicit expressions. The main use of Theorems 1 and 2 is the following consistency result.

THEOREM 3 (Consistency).    *If* (3.5) *holds, then the test which rejects the null hypothesis for large values of* $T_{\mathrm{w}}$ *is consistent against an alternative distribution* $\nu$ *if and only if*

$$E_{\nu} \left[ \frac{1}{f(x; \hat{\boldsymbol{\theta}}_{\nu})} \mathbf{t}(x) \right] \neq \mathbf{0}.$$

*In particular, the test is consistent against all alternatives if and only if* $a_k \neq 0$ *for all* $k$.

## 4. The rotation group $SO(3)$.

4.1. *Sobolev tests on $SO(3)$.* Two important Sobolev tests of uniformity on the rotation group $SO(3)$ are Downs' [4] generalization of the Rayleigh test and Prentice's [15] generalization of Giné's [7] $G_n$ test. See Section 13.2.2 of [14]. For a sample $\mathbf{X}_1, \dots, \mathbf{X}_n$ on $SO(3)$, these tests reject uniformity for large values of the Rayleigh statistic

$$T_{\mathrm{R}} = 3n \operatorname{tr}(\bar{\mathbf{X}}'\bar{\mathbf{X}}),$$

where

$$\bar{\mathbf{X}} = \frac{1}{n} \sum_{i=1}^{n} \mathbf{X}_i,$$

and the Giné statistic

$$T_{\mathrm{G}} = \frac{1}{n} \sum_{i=1}^{n} \sum_{j=1}^{n} \left( \frac{1}{2} - \frac{3\pi}{32} [\operatorname{tr}(\mathbf{I}_3 - \mathbf{X}_i'\mathbf{X}_j)]^{1/2} \right),$$

respectively. The corresponding goodness-of-fit tests reject the null hypothesis for large values of the weighted Rayleigh statistic

$$T_{\mathrm{wR}} = 3n \operatorname{tr}(\bar{\mathbf{X}}_{\mathrm{w}}' \bar{\mathbf{X}}_{\mathrm{w}}),$$

where

$$\bar{\mathbf{X}}_{\mathrm{w}} = \frac{1}{n} \sum_{i=1}^{n} \frac{1}{f(\mathbf{X}_i; \hat{\boldsymbol{\theta}})} \mathbf{X}_i,$$

and the weighted Giné statistic

$$T_{\mathrm{wG}} = \frac{1}{n} \sum_{i=1}^{n} \sum_{j=1}^{n} \frac{1}{f(\mathbf{X}_i; \hat{\boldsymbol{\theta}}) f(\mathbf{X}_j; \hat{\boldsymbol{\theta}})} \left( \frac{1}{2} - \frac{3\pi}{32} [\operatorname{tr}(\mathbf{I}_3 - \mathbf{X}_i'\mathbf{X}_j)]^{1/2} \right),$$

respectively. For $T_{\mathrm{R}}$ and $T_{\mathrm{wR}}$, $a_k = 0$ for $k \geq 2$; for $T_{\mathrm{G}}$ and $T_{\mathrm{wG}}$, all the $a_k$ are nonzero ([15], pages 173–174). It follows from Theorem 3 that the goodness-of-fit test based on $T_{\mathrm{wR}}$ is consistent only against alternatives $\nu$ with $E_{\nu}[\mathbf{X}] \neq \mathbf{0}$, whereas the test based on $T_{\mathrm{wG}}$ is consistent against all alternatives.

4.2. *A numerical example.* The set of vectorcardiogram data described in [5] is a classic data set on $SO(3)$. The portion of this data set given by the orientations of vectorcardiograms obtained using the Frank lead system from boys aged 2–10 gives 28 observations on $SO(3)$. For these 28 observations, $T_{\mathrm{R}} = 209$, so that comparison with the large-sample limiting $\chi_9^2$ distribution (which is appropriate for $n \geq 18$ by Table 1 of [8]) indicates very clearly that uniformity should be rejected.

The eigenvalues of $\bar{\mathbf{X}}$ are $0.957, 0.888$ and $0.883$, suggesting that it is appropriate to fit a matrix Fisher distribution with canonical parameter matrix of the form $\kappa\mathbf{U}$, where $\kappa > 0$ and $\mathbf{U} \in SO(3)$, that is, the probability density function is

$$f(\mathbf{X}; \mathbf{U}, \kappa) = M(\tfrac{1}{2}, 2, 4\kappa)^{-1} e^{\kappa} \exp\{\kappa \operatorname{tr}(\mathbf{U}'\mathbf{X})\},$$

where $M(1/2, 2, \cdot)$ is a Kummer function. (See [4, 12], or Section 13.2.3 of [14].) The maximum likelihood estimates of $\kappa$ and $\mathbf{U}$ are $\hat{\kappa} = 5.63$ and

$$\hat{\mathbf{U}} = \begin{pmatrix} 0.583 & 0.629 & 0.514 \\ 0.660 & -0.736 & 0.151 \\ 0.473 & 0.252 & -0.844 \end{pmatrix}.$$

The $p$-values (based on 1000 simulations) of the goodness-of-fit tests are 0.169 for the weighted Rayleigh test and 0.126 for the weighted Giné test, each indicating clearly that the fit is acceptable.

**Acknowledgments.** I am grateful to Professor T. D. Downs for giving me access to the vectorcardiogram data and to a referee for the suggestion of allowing the mapping $\mathbf{t}$ to depend on the sample size.

SCHOOL OF MATHEMATICS AND STATISTICS
UNIVERSITY OF ST. ANDREWS
NORTH HAUGH
ST. ANDREWS KY16 9SS
UK
E-MAIL: pej@st-andrews.ac.uk